\documentclass[12pt]{amsart}
\usepackage{amsmath,amssymb,amsthm,verbatim,url,geometry}
\usepackage{rotating}

\newcommand{\nc}{\newcommand}
\nc{\transpose}[1]{#1^\mathrm{t}}
\nc\T{\rule{0pt}{3.1ex}} \nc\B{\rule[-1.7ex]{0pt}{0pt}} 
\nc{\V}{{\{0,1\}^*}}
\nc{\Impl}{\Rightarrow}
\nc{\Rys}{\op{Rys}}
\nc{\rev}{\mathrm{r}}
\newtheorem{thm}{Theorem}\nc{\bthm}{\begin{thm}} \nc{\ethm}{\end{thm}}
\newtheorem{prop}[thm]{Proposition}\nc{\bprp}{\begin{prop}} \nc{\eprp}{\end{prop}}
\newtheorem{fact}[thm]{Fact}\nc{\bfct}{\begin{fact}} \nc{\efct}{\end{fact}}
\newtheorem{prob}[thm]{Problem}\nc{\bprb}{\begin{prob}} \nc{\eprb}{\end{prob}}
\newtheorem{lem}[thm]{Lemma}\nc{\blem}{\begin{lem}} \nc{\elem}{\end{lem}}
\newtheorem{claim}[thm]{Claim}\nc{\bclm}{\begin{claim}} \nc{\eclm}{\end{claim}}
\newtheorem{cor}[thm]{Corollary}\nc{\bcor}{\begin{cor}} \nc{\ecor}{\end{cor}}
\newtheorem{conj}[thm]{Conjecture}\nc{\bcnj}{\begin{conj}} \nc{\ecnj}{\end{conj}}
\theoremstyle{definition}\newtheorem{defn}[thm]{Definition}\nc{\bdfn}{\begin{defn}} \nc{\edfn}{\end{defn}}
\theoremstyle{remark}
\newtheorem{rem}[thm]{Remark}\nc{\brem}{\begin{rem}}\nc{\erem}{\end{rem}}
\newtheorem{cnv}[thm]{Convention}\nc{\bcnv}{\begin{cnv}} \nc{\ecnv}{\end{cnv}}
\newtheorem{exam}[thm]{Example}\nc{\bexm}{\begin{exam}}\nc{\eexm}{\end{exam}}
\newtheorem{alg}[thm]{Algorithm}\nc{\balg}{\begin{alg}}\nc{\ealg}{\end{alg}}
\nc{\bpf}{\begin{proof}}\nc{\epf}{\end{proof}}\nc{\be}{\begin{enumerate}}\nc{\ee}{\end{enumerate}}
\nc{\bi}{\begin{itemize}}\nc{\itm}{\item}\nc{\ei}{\end{itemize}}
\nc{\Cayley}{\op{Cayley}}
\nc{\x}{\times}
\nc{\inv}{^{-1}}
\nc{\cT}{\mathcal{T}}
\nc{\cK}{\mathcal{K}}
\nc{\cD}{\mathcal{D}}
\nc{\tr}{\op{tr}}
\nc{\sm}{\setminus}
\nc{\sub}{\subseteq}
\nc{\set}[2]{\left\{#1\;:\;#2\right\}}
\nc{\ceq}{\stackrel{\cdot}{=}}
\nc{\my}[1]{\textcolor{red}{\sf #1}}
\nc{\mx}[1]{\begin{pmatrix}#1\end{pmatrix}}
\nc{\op}[1]{\operatorname{#1}}
\nc{\SL}[1]{\SLtwo({#1})}
\nc{\bbF}{\mathbb{F}}
\nc{\SLtwo}{\op{SL}_2}
\nc{\SLq}{$\mathrm{{SL}}_2(\bbF_q)$}
\nc{\SLx}{\mathrm{{SL}}_2(\bbF_q[x])}
\nc{\Mat}[2]{\mathrm{{M}}_{#1}({#2})}
\nc{\GL}{\mathrm{GL}}
\nc{\bM}{\begin{smallmatrix}}
\nc{\eM}{\end{smallmatrix}}
\nc{\ed}{\end{document}}

\title{On Rystov's generalization of the \v{C}ern\'y Conjecture}

\author{Noam Lifshitz}
\address{Department of Mathematics, Bar Ilan University, 5290002 Ramat Gan, Israel}
\email{noamlifshitz@gmail.com}

\author{Ciaran Mullan}
\address{Technische Universit\"at Darmstadt, Fachbereich Informatik, Kryptographie und Computeralgebra,
Hochschulstra\ss{}e 10, 64289 Darmstadt, Germany}

\author{Boaz Tsaban}
\address{Department of Mathematics, Bar Ilan University, 5290002 Ramat Gan, Israel}
\email{tsaban@math.biu.ac.il}
\urladdr{http://www.cs.biu.ac.il/\~{}tsaban}

\begin{document}

\begin{abstract}
We resolve a conjecture of Rystov concerning products of matrices,
that generalizes the \v{C}ern\'y Conjecture. This is a preliminary announcement. Later versions
will include additional results and details.
\end{abstract}

\maketitle

The following conjecture is a generalization of the celebrated, and still open,
\v{C}ern\'y Conjecture.
By \emph{product in $A_1,\dots,A_k$ of length $l$} we mean a product
$$A_{i_1}A_{i_2}\cdots A_{i_l}$$
with $1\le i_1,\dots,i_l\le k$ not necessarily distinct.

\bcnj[Rystov]
Let $A_1,\dots,A_k\in\Mat{n}{\bbF}$, $k$ arbitrary. If the semigroup generated by $A_1,\dots,A_k$
is finite and contains the zero matrix $O$, then there is a product in $A_1,\dots,A_k$ of length at most $n^2$
that is equal to $O$.
\ecnj

Over finite fields, the condition that the generated semigroup is finite is satisfied automatically, and the conjecture
asserts that if any product in matrices $A_1,\dots,A_k$ equals $O$, then there is one of length at most $n^2$.
As there are no zero divisors in a field, the conjecture is true when $n=1$.
We prove that this conjecture fails for all $n>1$.
In the language of semigroup theory, the following lemma is equivalent to the
assertion that the semigroup of all singular matrices in $\Mat{2}{\bbF}$ is categorical at $0$. 

\blem[folklore]\label{lem:abc}
Let $A,B,C\in\Mat{2}{\bbF}$, $B$ singular. If $ABC=O$, then $AB=O$ or $BC=O$.
\elem
\bpf
If $C$ is invertible we are done. Thus, assume that $C$ is singular.
If $C=O$ we are done, so assume $C$ is nonzero.
Then $0$ is an eigenvalue of $C$, and the characteristic polynomial of $C$ is $x(x-\gamma)=0$.
If $\gamma=0$ then the Jordan form of $C$ is 
$$\mx{0 & 1\\ 0 & 0}.$$
If $\gamma\neq 0$, then $C$ is conjugate to a nonzero diagonal matrix of the form
$$\mx{* & 0\\ 0 & 0}.$$
Since scalar multiplication does not affect our problem, we may assume in this case that 
$C$ is conjugate to 
$$\mx{1 & 0\\ 0 & 0}.$$
Conjugate $A,B,C$ by the same matrix (and multiply $C$ by a nonzero scalar, if needed), 
such that $C$ obtains one of the two above-mentioned forms.

As $B$ is singular, there is $v$ such that the columns of $B$ are $(v,\beta v)$
or $(\vec{0},v)$. If $B=(v,\beta v)$, then
$$O=ABC = A(v,\beta v)C = (Av,\beta Av)C = (Av,\vec{0})\mbox{ or } (\vec{0},Av),$$
depending on the form of $C$. Then $Av=\vec{0}$, and therefore $AB=(Av,\beta Av)=O$.
And if $B=(\vec{0},v)$, then
$$BC = (\vec{0},v)\mx{1 & 0\\ 0 & 0} = O.\qedhere$$
\epf

Lemma \ref{lem:abc} has the following consequence.
Assume, in general, that 
$$A_1\cdots A_k=O$$
in $\Mat{2}{\bbF}$.
By cancelling off the regular matrices on the edges,
We may then assume that $A_1$ and $A_k$ are singular. 
Next, if any $A_i$ is singular, $1<i<k$, we can apply
Lemma \ref{lem:abc} to 
$$(A_1\cdots A_{i-1})A_i(A_{i+1}\cdots A_{k-1}A_k)$$
and conclude that 
$$(A_1\cdots A_{i-1})A_i=O\mbox{ or } A_i(A_{i+1}\cdots A_{k-1}A_k)=O.$$
We can continue this procedure until we arrive at a word of the form
$$A_1\cdots A_k=O$$
where $A_1$ and $A_k$ are singular, and all other matrices are invertible.

\bcor\label{cor:abba}
Let $A,B\in\Mat{2}{\bbF}$ with $B$ invertible. If any product in $A,B$ equals $O$, then
there is a unique shortest product in $A,B$ that equals $O$. The shortest product is of the
form $AB^mA=O$. \qed
\ecor

We are now ready to prove the main result.

\bthm\label{thm:main}
Let $n\ge 2$. For each $N$, there is a finite field $\bbF$ and matrices $A,B\in\Mat{n}{\bbF}$ such that
the minimal length of a product in $A,B$ that equals $O$ exists, and its length is greater than $N$.
\ethm
\bpf
If $A,B\in\Mat{2}{\bbF}$ exemplify the assertion for $n=2$, then for
every larger $n$, the block matrices
$$\mx{A & O\\O & O},
\mx{B & O\\O & O}\in\Mat{n}{\bbF}$$
exemplify the same assertion. Thus, we may assume that $n=2$.

Take 
$$A=\mx{1 & 0\\0 & 0}, B=\mx{1 & 1\\1 & 0},$$
and note that $B$ is invertible and $A$ is idempotent.
As 
$$B\mx{y\\x}=\mx{x+y\\y},$$
we have that
$$B^k\mx{1\\0}=\mx{F_{k+1}\\F_k}$$
for all $k=0,1,\dots$, where $F_0=0, F_1=1$, and $F_{k+1}=F_k+F_{k-1}$
for all $k>1$, that is, $F_k$ is the $k$-th element of the Fibonacci sequence.
Thus, 
$$AB^kA=\mx{1 & 0\\0 & 0}B^k\mx{1 & 0\\0 & 0}=\mx{1 & 0\\0 & 0}\mx{F_{k+1} &0\\ F_k & 0}=F_{k+1}$$
for all $k$.
Fix a prime number $p\ge F_{N+1}$. 
Let $k$ be minimal such that $F_{k+1}=0\pmod p$.
Necessarily, $k\ge N$.
Then, in $\Mat{2}{\bbF_p}$,
$$AB^kA=O,$$
and $AB^mA\neq O$ for all $0\le m<k$.
By Corollary \ref{cor:abba}, the length of the shortest $O$ product is $k+2\ge N+2>N$.
\epf

\subsection*{Acknowledgments}
We learned of the Rystov Conjecture from Benjamin Steinberg, who presented it in an
open problem session during a recent conference in honor of Stuart Margolis.
We thank Ben and Stuart for stimulating discussions on this problem.

\end{document}